\documentclass[12pt]{amsart}

\usepackage{amssymb,latexsym}

\usepackage{enumerate}

\makeatletter

\@namedef{subjclassname@2010}{%

  \textup{2010} Mathematics Subject Classification}

\makeatother

\newtheorem{theorem}{Theorem}[section]

\newtheorem{corollary}[theorem]{Corollary}

\newtheorem{lemma}[theorem]{Lemma}

\theoremstyle{definition}

\numberwithin{equation}{section}

\newcommand{\subgp}[1]{\langle{#1}\rangle}

\newcommand{\qfd}{\hfill $\Box$}

\begin{document}


\baselineskip=17pt




\date{}

\title[ Minkowski Product Size]{On  Minkowski product size: \\The Vosper's property }
\author[ Y. O. Hamidoune]{ Yahya Ould Hamidoune}
\address{UPMC Univ Paris 06,
 E. Combinatoire, \\Case 189, 4 Place Jussieu,
75005 Paris, France}
\email{hamidoune@math.jussieu.fr}

\maketitle


\begin{abstract}
A subset $S$ of a group $G$ is said to be a Vosper's subset if
$|A\cup AS|\ge \min (|G|-1,|A|+|S|),$ for any subset $A$ of $G$ with $|A|\ge 2.$ In the present work, we describe  Vosper's subsets. Assuming  that $S$ is not a progression
and that  $|S^{-1} S|, |S S^{-1}| <2 |S|,|G'|-1,$ we show that
   there exist  an element $a\in S,$ and a non-null subgroup $H$ of $G'$ such that
 either
 $S^{-1}HS =S^{-1}S \cup a^{-1}Ha$ 
 or
$SHS^{-1} =SS^{-1}\cup aHa^{-1},$ where $G'$ is the subgroup generated by $S^{-1}S.$

\end{abstract}



\maketitle

\section{Introduction}
Let  $ A,B$ be
subsets of a group $ G $. The {\em Minkowski product} of $ A$ with $B$ is defined as
$$AB=\{xy \ : \ x\in A\  \mbox{and}\ y\in
  B\}.$$

  Kneser's Theorem \cite{knesrcomp} states that  $AB$ is a periodic set if $|AB|\le |A|+|B|-2$ and if  $G$ is abelian. Diderrich \cite{diderrich} obtained the same conclusion assuming only that the elements of $B$ commute. As mentioned in \cite{hiso2007}, the last result follows from Kneser's Theorem. In \cite{olsonsdif}, Olson constructed subsets $A$ and $B$ of some non-abelian group with $|AB|\le |A|+|B|-2$ such that for every  non-null group $H,$

  $AB\neq AHB$,  $AB\neq HAB$ and $AB\neq ABH.$

  The special cases $B=A$ and $B=A^{-1}$ received also some attention.  In \cite{freiman},
 Freiman described subsets $A$ with  $|A^2|< \frac{1+\sqrt{5}}{2}|A|$ or $|A^{-1}A|< \frac{1+\sqrt{5}}{2}|A|.$
 A transparent  exposition of Freiman results is contained in Husbands dissertation \cite{husb}.

  Tao proposed in \cite{t2} a short proof of Freiman's result, suggesting that threshold should be ${2}.$ In \cite{hkft},  we obtained a Kneser type result, asserting  that there exists is a non-null subgroup $H$ of $G$ such that
 $A^{-1}HA =A^{-1}A$ or $AHA^{-1} =AA^{-1},$ if  $|A^{-1}A |< 2|A|-1.$

As mentioned by Tao, in \cite{t2}, the relations $|A^{-1}A |< 2|A|$ and $|AA^{-1} |< 2|A|$ imply no kind of periodicity, since they are satisfied by left-progressions in a torsion free  group.
  The methods  used in \cite{hkft} are not enough precise to give an inverse theorem for  $|A^{-1}A |< 2|A|.$

We need to develop the isoperimetric approach in the non-abelian
 case, continuing the work done in the finite case in \cite{hast}. The present work  generalizes  the results obtained  in the abelian case \cite{hejc3} and  the results obtained in the non-abelian  finite case in \cite{hast,hactaa}. The  two papers \cite{hejc3,hast} use the obsolete language of super-atoms. We use here the more general and more precise language of $k$-atoms introduced in \cite{halgebra}. Instead of restricting ourselves to the case of a Minkowski product,  we develop the approach for an arbitrary relation. The information
 on Minkowski product will follow, once we restrict ourselves to Cayley relations $x^{-1}y\in A.$ In almost all cases, the results obtained in the special case of Cayley relations hold for relations having a transitive group of automorphisms.

 Among other tools, the isoperimetric approach,  was used by Serra-Zemor \cite{sz} and by Vu-Wood \cite{vw} to replace the classical rectification. It was also used by the author \cite{hkemp+1} to propose a geometric approach to the classical Kemperman Theory \cite{kempacta}, leading to simplifications and generalizations.

Let $\Gamma=(V,E)$ be a reflexive relation. The board of a subset $X$ is $\Gamma (X)\setminus X.$  Put ${\mathbf F}_k=\{X\subset V :|X|\ge k \ \text{and} \ |V\setminus \Gamma (X)|\ge k\}.$
The $kth$-connectivity $\kappa _k$ is the minimal cardinality of the boards of the members of
${\mathbf F}_k.$ A member of ${\mathbf F}_k$ achieving this minimum will be called a $k$-fragment. A $k$-fragment with minimal (resp. maximal when $V$ is finite) cardinality
will be called a $k$-atom (resp. $k$-super-fragment). The relation will be called $k$-faithful if $|A|\le |V\setminus \Gamma (A)|,$
 where $A$ is a $k$-atom. A relation $\Gamma$ will be called a Cauchy relation if $\kappa _1 \ge |\Gamma (v)|-1,$ for some $v\in V.$
  A relation with $\kappa _2>\kappa _1$ will be called a Vosper's relation. In this language, the Cauchy-Davenport Theorem \cite{cauchy,davenport} states that Cayley relations on groups with a prime order are Cauchy relations. Vosper's Theorem \cite{vosper1} states that, for $|G|$  a prime,  the Cayley relation $x^{-1}y\in A$  is a Vosper's relation, if $A$ is not an arithmetic progression.

 Our main problem is to describe the Vosper's Cayley relations.The organization of the paper is the following:

 Section 2 contains some terminology. The basic notions are presented in section 3. In Section 4,  we prove  that the intersection of distinct $k$-atoms of a $k$-faithful relation has cardinality less than $k.$  We show also in this section that the intersection of distinct $k$-super-fragments  has cardinality less than $k,$
  when the reverse relation is non-$k$-faithful.
  In section 5, we obtain more precise intersection properties for non Vosper's relations.
  In section 6, we investigate the intersection of three $2$-atoms.
 In section 7, we apply the last result to describe Vosper's relations with a transitive group of automorphisms. In section 8, we show that one of the two Cayley relations  $x^{-1}y\in A$ and $x^{-1}y\in A^{-1}$ has  $2$-atom of the form $H\cup Ha,$ where $H$ is a subgroup and $a$ is an element of $G.$ As an application, we obtain
 in section 9, we obtain the following result:

   \begin{theorem}\label{kneser}
Let $ A$ be a subset of group $G_0$ and let $G$ be the subgroup generated by
${A^{-1}A}.$  If  $|A^{-1} A|, |A A^{-1}| <2 |A|,$ then
one of the following holds:
\begin{itemize}
  \item[(i)] There is an $u\in G$ with $u^2=1$ such that either $A A^{-1}=G\setminus \{u\}$
  or $A^{-1}A=G\setminus \{u\},$
  \item[(ii)]  $A$ is a progression,
  \item[(iii)]   there exists is a non-null subgroup $H$ of $G$ such that
  $A^{-1}HA =A^{-1}A\cup a^{-1}Ha,$ for some $a\in  A,$
 \item[(iv)]   there exists is a non-null subgroup $H$ of $G$ such that
$AHA^{-1} =AA^{-1}\cup aHa^{-1},$ for some $a\in A.$
\end{itemize}

\end{theorem}

\section{Some Terminology}

 An ordered pair $\Gamma =(V,E),$ where $V$ is a set and  $E\subset V\times V,$ will be called a {\em graph}
or a {\em relation} on  $V.$Let $\Gamma =(V,E)$ be a  graph and let $X\subset V.$ The {\em reverse} graph of $\Gamma$ is
the graph $\Gamma ^- =(V,E^-),$ where $E^-=\{(x,y): (y,x)\in E\}.$
The graph $\Gamma$ will be called {\em locally-finite}  if for all $x\in V,$ $|\Gamma (x)|$ and $|\Gamma ^-(x)|$ are finite.
The graph $\Gamma$ is said to be {\em $r$-regular} if $|\Gamma (x)|=r,$ for every $x\in V.$
The graph $\Gamma$ is said to be {\em $r$-reverse-regular} if $|\Gamma ^-(x)|=r,$ for every $x\in V.$
The graph $\Gamma$ is said to be {\em $r$-bi-regular} if it is {\em $r$-regular} and {\em $r$-reverse-regular}.

 \begin{itemize}
  \item The minimal degree of $\Gamma$ is defined as $\delta ({\Gamma})= \min \{ |\Gamma (x)| : x\in V\}.$
  \item We write $\delta _{\Gamma ^-}= \delta _{-} (\Gamma).$
  \item The board of  $X$  is defined as $\partial _{\Gamma}(X)= \Gamma (X)\setminus X.$
  \item The exterior of  $X$  is defined as
 $\nabla _{\Gamma}(X)= V\setminus \Gamma (X).$
  \item We shall write $\partial ^{-}_{\Gamma } =\partial _{\Gamma ^{-}}.$ This subset will be called the
  {\em reverse-board} of  $X.$
  \item  We shall write $\nabla ^{-}_{\Gamma } =\nabla _{\Gamma ^{-}}.$
\end{itemize}

When the context is clear, the reference to $\Gamma$ will be omitted.

\section{Basic notions}

In this section, we define the concepts of $kth$-connectivity,
$k$-fragment and $k$-atom and prove some elementary properties of
these notions.

A graph $\Gamma $ will be called $k$-separable if there is a finite subset  $X\subset V,$ with
  $k\le |X|<\infty$  { and} $k\le |V\setminus  \Gamma(X)|.$
 The {\em $kth$-connectivity} of a $k$-separable graph
 $\Gamma$ (called  {\em $kth$-isoperimetric number} in
\cite{halgebra})
 is defined  as
\begin{equation}
\kappa _k (\Gamma )=\min  \{|\partial (X)|\   :  \ k\le |X|<\infty \  \text{ and} \ k\le |V\setminus  \Gamma(X)|
\}.
\label{eqcon}
\end{equation}

 A finite subset $X$ of $V$ such that $k\le |X|<\infty,$   $k\le |V\setminus  \Gamma(X)|$ and $|\partial (X)|=\kappa
_k(\Gamma)$ is called a {\em $k$-fragment} of $\Gamma$. A
$k$-fragment with minimum cardinality is called a {\em $k$-atom}.

These notions were introduced in \cite{halgebra}. Let us now introduce more notions.
A  subset $X$ of $V$ will be called a {\em $k$-semi-fragment} of $\Gamma$ if
either $X$ is a $k$-fragment or
 $\nabla (X)$ is a reverse $k$-fragment.
A
$k$-fragment of a finite graph having a maximal cardinality will be called a {\em $k$-super-fragment}.
 The graph $\Gamma$ will be called $k$-{\em faithful} if
$|A|\le |\nabla (A)|,$ where $A$ is a k-atom.

A $k$-semi-fragment  of $\Gamma ^{-}$ will be  called
a {\em reverse}-$k$-semi-fragment of $\Gamma$.
A $k$-fragment  of $\Gamma ^{-}$ will be  called a {\em reverse}
$k$-fragment of $\Gamma$. We shall write $\kappa_{-k}(\Gamma ) =\kappa_k(\Gamma ^{-}).$ The reference
to $\Gamma$ could be implicit.

Recall that $\kappa _k (\Gamma)$ is the maximal integer $j$
such that for every finite subset $X\subset V$  with $|X|\geq k$,
$|\Gamma (X)|\geq \min \Big(|V|-k+1,|X|+j\Big).$

The following lemma is immediate from the definitions:

\begin{lemma} \label{degenerate}{ Let $k\ge 2$ be an integer. A  reflexive
  locally finite   $k$-separable graph  $\Gamma =(V,E)$ is  is a $k-1$-separable graph, and moreover
  $\kappa _{k-1}\le \kappa _k.$ If $\kappa _{k-1}= \kappa _k,$  then  $${\mathbf F}_k =\{F\in  {\mathbf F}_{k-1} :
  k\le \min (|F|,|\nabla (F)|)\}.$$

  }
  \end{lemma}


The next lemma contains useful duality relations:

\begin{lemma} \label{negative}{Let  $X$ and $Y$ be $k$-fragments of a reflexive
  locally finite   $k$-separable graph $\Gamma =(V,E)$.  Then
 \begin{align}
\partial ^{-}  (\nabla (X))&=
\partial (X),\label{eqduall}\\
\nabla^- (\nabla (X))&={X}, \label{eqdualf}\\
X\subset Y & \text{if and only if} \ \nabla (Y)\subset \nabla (X). \label{nomon}
 \end{align}

 In particular,  $\nabla (X)$ is a reverse-$k$-semi-fragment.

}\end{lemma}
\proof Clearly, $ \partial (X) \subset \partial ^{-} (\nabla (X))$

We must have $ \partial (X) = \partial ^{-} (\nabla (X))$, since
otherwise there is a $y\in \partial ^{-} (\nabla (X))
\setminus \partial (X).$ It follows that
$|\partial (X\cup \{y\})|\le |\partial (X)|-1$, contradicting the
definition of $\kappa _k.$ This proves (\ref{eqduall}). In particular, $\nabla (X)$ is a reverse-$k$-semi-fragment.

Thus $\Gamma ^{-}
(\nabla (X))=\nabla (X)\cup
\partial ^{-}(\nabla (X)) =\nabla (X)\cup
\partial (X)=V\setminus X.$
Thus (\ref{eqdualf}) holds.
Clearly, (\ref{nomon}) is a direct consequence of (\ref{eqdualf}).\qfd

Let us  define two important notions:

Let $\Gamma =(V,E)$ be a  reflexive graph. We shall say that
$\Gamma$ is a {\em Cauchy graph}  if $\Gamma$ is  non-$1$-separable or if $\Gamma$ has a $1$-atom $A$
with $|A|=1$ or $|\nabla (A)|=1.$  We shall say that
$\Gamma$ is a {\em reverse-Cauchy graph} if $\Gamma ^-$ is a Cauchy graph.

Clearly, $\Gamma$ is a Cauchy graph if and only if for every
$X\subset V$ with $|X|\ge 1$,
 $$|\Gamma(X)|\ge \min \Big(|V|, |X|+\delta-1\Big).$$

 We shall say that
$\Gamma$ is  {\em degenerate} if $\Gamma$ is  $2$-separable and $\kappa _2= \kappa _1.$
We shall say that
$\Gamma$ is  {\em reverse-degenerate} if $\Gamma ^-$ is degenerate.


\begin{lemma} \label{finiteg}{Let $\Gamma =(V,E)$ be a reflexive
finite   $k$-separable graph and let $X$  be a subset of $V.$
Then \begin{equation}
\kappa _k= \kappa _{-k}.\label{kk-}\end{equation}  Moreover,
\begin{itemize}
  \item[(i)]   $X$ {is a } $k$-{fragment } {if and only if } $\nabla(X)$ {is a } $k$-{reverse-fragment,}
  \item[(ii)] $X$ {is a } $k$-{super-fragment } {if and only if } $\nabla(X)$ {is a } $k$-{reverse-atom,}
  \item[(iii)] $\Gamma$ is a Cauchy graph if and only if it is a reverse-Cauchy graph,
   \item[(iv)] $\Gamma$ is degenerate if and only if it is reverse-degenerate.
\end{itemize}
}\end{lemma}
\proof
Observe that a finite graph is $k$-separable if and only if its reverse
is $k$-separable.
Take a
$k$-fragment $X$  of $\Gamma$. We have clearly $\partial _{-}
(\nabla (X))\subset
\partial (X)$.
Therefore $$\kappa _k(\Gamma )\ge |\partial (X)|\ge |\partial ^{-}
(\nabla (X))|\ge \kappa _{-k}.$$ The reverse inequality of (\ref{kk-}) follows similarly or by duality.

Suppose that $X$ is a $k$-fragment.
By (\ref{eqduall}) and (\ref{kk-}), $|\partial _{-}
(\nabla (X))|=
|\partial (X)|=\kappa _{k}=\kappa _{-k},$ and hence $\nabla (X)$ is a revere $k$-fragment. The other implication of
(i) follows similarly.
Suppose now that $X$ is a $k$-super-fragment. By (i), $\nabla (X)$ is a reverse-$k$-fragment. Take a reverse-$k$-atom $N.$
Now $\nabla ^- (N)$ is a $k$-fragment by (i). Thus $|\nabla ^- (N)|\le |X|,$ and hence using (\ref{eqdualf}), $|\nabla  (X)|\le |N|.$
Thus, $\nabla (X)$ is a reverse-$k$-atom. The other implication of
(ii) follows similarly.
Now (iii) and (iv) follow directly from the definitions and (\ref{kk-}).\qfd

Recall the following easy fact:
\begin{lemma}\label{katomdegree}\cite{halgebra}
Let $\Gamma =(V,E)$ be a  locally-finite $k$-separable graph and let $A$ be a $k$-atom with $|A|>k.$
 Then $\Gamma ^-(x)\cap A\neq \{x\},$ for every $x\in A.$ 
\end{lemma}

\proof We can not have $\Gamma ^-(x)\cap A= \{x\},$ otherwise $ A\setminus \{x\}$
would be a $k$-fragment.\qfd

\section{  Geometric properties of fragments}

The next result generalizes  results obtained in \cite{hast,halgebra,hiso2007}:

\begin{theorem}
Let  $X$ be a $k$-fragment of
 a  reflexive locally finite
$k$-separable graph   $\Gamma =(V,E)$
and let $Y$ be a  $k$-semi-fragment.

\begin{itemize}
  \item[(i)] If $|X\cap Y|\ge k,$ then $|\nabla (Y)\cap \partial (X)|\le |X\cap \partial (Y)|,$
  \item[(ii)]If $|X\cap Y|\ge k$ and $|\nabla (X)\cap \nabla (Y)|\ge k,$ then $X\cap Y$ is a $k$-fragment
   \item[(iii)] If $|X\cap Y|\ge k$ and $|X|\le | \nabla (Y)|,$ then $X\cap Y$ is a $k$-fragment
\end{itemize}
\label{lem1977} \end{theorem}

\proof
{\begin{center}
\begin{tabular}{|c||c|c|c|c|}
\hline
$\cap $&  $Y$ & $\partial (Y) $ &  $\nabla (Y)$ \\
\hline
\hline
$X$&  $R_{11}$ & $R_{12}$ & $R_{13}$ \\
\hline
$\partial (X) $& $R_{21}$& $R_{22}$ & $R_{23}$ \\
\hline
$\nabla (X)$ & $R_{31}$  & $R_{32}$ & $R_{33}$ \\
\hline
\end{tabular}

\end{center}}

Assume that $|X\cap Y|\ge k.$
 By the definition of $\kappa_k,$
\begin{align*}  |R_{21}|+|R_{22}|+ |R_{23}| &= \kappa_k
\\&\le |\partial (X\cap Y)|\\
&=  |R_{12}|+|R_{22}|+ |R_{21}|,
\end{align*}
and hence $|\nabla (Y)\cap \partial (X)|=|R_{23}|\le
|R_{12}|=|X\cap \partial (Y)|,$ showing (i).

We shall prove that
\begin{equation} \text{If}\ |\nabla (X)\cap \nabla (Y)|\ge k,\text{then}\ |\nabla (Y)\cap \partial (X)|\ge |X\cap \partial (Y)|,\label{iia}
\end{equation}

Assuming that $ |\nabla (X)\cap \nabla (Y)|\ge k,$

{\bf Case} 1: $Y$ is finite.
\begin{align*}  |R_{12}|+|R_{22}|+ |R_{32}| &=  \kappa_{k}
\\&\le |\partial (X\cup Y)|\\
&\le  |R_{22}|+|R_{23}|+ |R_{32}|,
\end{align*}
and hence $|R_{12}|\le |R_{23}|,$ showing (\ref{iia}) in this case.

{\bf Case} 2: $Y$ is infinite.
\begin{align*}  |R_{12}|+|R_{22}|+ |R_{32}| &=  \kappa_{-k}
\\&\le |\partial ^- (R_{33})|\\
&\le  |R_{22}|+|R_{23}|+ |R_{32}|,
\end{align*}
and hence $|R_{12}|\le |R_{23}|,$ showing (\ref{iia}) in this case.

Assume now that $|X\cap Y|\ge k$ and $|\nabla (X)\cap \nabla (Y)|\ge k.$ By (i) and (\ref{iia}),
we have $|\nabla (Y)\cap \partial (X)|= |X\cap \partial (Y)|.$ It follows that
   $$\kappa _k\le |\partial (X\cap Y)|\le  |R_{12}|+|R_{22}|+ |R_{21}|\le  |R_{12}|\le |R_{23}|+|R_{22}|+ |R_{21}|=\kappa _k,$$
showing that $X\cap Y$ is a $k$-fragment. 

Assume now that  $|X\cap Y|\ge k$ and $|X|\le | \nabla (Y)|.
$ By (i), $|R_{12}|\ge |R_{23}|.$
Clearly, $$|R_{13}|+|R_{23}|+R_{33}= |\nabla (Y)|\ge |X|=|R_{11}|+|R_{12}|+R_{13}.$$

Therefore $ |\nabla (X)\cap \nabla (Y)|\ge |X\cap Y|\ge k.$ Now, (iii) by applying (ii).\qfd

We shall now investigate  the  super-fragments behavior when the atoms are too big.
Let us  mention two easy facts:

\begin{lemma}\label{fourf}
A $k$-separable graph $\Gamma=(V,E)$ is either $k$-faithful or reverse $k$-faithful. Moreover infinite graphs are $k$-faithful.
 \end{lemma}
\proof Assume that $\Gamma=(V,E)$ is non-$k$-faithful. Then $V$ is clearly finite. Let $A'$ be a reverse $k$-atom.
By Lemma \ref{finiteg}, $\nabla ^- (A')$ is a $k$-fragment and $\nabla (A)$ is a reverse $k$-fragment.
By (\ref{nomon}), we have $|\nabla ^- (A')|\ge |A|> |\nabla (A)|\ge |A'|.$
 In particular $\Gamma$ is reverse-faithful.\qfd

\begin{theorem}

Let $\Gamma =(V,E)$ be a  reflexive  finite
$k$-separable graph such that  $\Gamma^- $ is a non-$k$-faithful graph. Then

 \begin{itemize}
   \item[(i)] the intersection of two distinct $k$-super-fragments has a cardinality less than  $k.$
   \item[(ii)] Moreover,  if $k\ge 2$ and $\kappa _k=\kappa _{k-1},$
then the intersection of two distinct $k$-super-fragments has a cardinality
 less than  $k-1.$
  \end{itemize}

\label{antiatom} \end{theorem}
\proof
Let $X$ and $Y$ be $k$-super-fragments of
$\Gamma.$
 By Lemma \ref{finiteg},  $\nabla (X)$ and  $\nabla (Y)$ are  reverse $k$-atoms  of $\Gamma.$
Since $\Gamma ^-$ is non-$k$-faithful, we have by (\ref{eqdualf}), $|\nabla (X)|>|\nabla ^- (\nabla (X))|=|X|.$

Suppose that  $|X\cap Y|\ge k.$
By Theorem \ref{lem1977},(iii), $X\cap Y$ is a $k$-fragment and hence $X=Y,$ a contradiction.

Assume now that $\kappa _k=\kappa _{k-1}$ and that $|X\cap Y|\ge k-1.$
By Lemma \ref{degenerate}, $X$ and $Y$ are $k-1$-fragments.

 By Theorem \ref{lem1977},(i), $|\nabla (Y)\cap \partial (X)|\ge |X\cap \partial (Y)|.$
Thus,
\begin{align*}
k-1\le |X\cap Y|&= |X| -|X\cap \partial  (Y)|-
|X\cap \nabla (Y)|\\&\le |\nabla (Y)|-1-|\nabla (Y) \cap  \partial  (X)|-|X\cap \nabla (Y)|\\
&=|\nabla (X)\cap \nabla (Y)|-1.
\end{align*}
 By Theorem \ref{lem1977},(ii), applied to $\Gamma ^-,$  $\nabla (X)\cap \nabla (Y)$ is a $k-1$-reverse-fragment.
By Lemma \ref{degenerate}, $\nabla (X)\cap \nabla (Y)$ is a $k$-reverse-fragment. Thus, $\nabla (X)=\nabla (Y)$, and hence $X=Y,$
a contradiction.\qfd

\section{Degenerate graphs}

The next consequence of Theorem \ref{lem1977} will be a main tool:

\begin{theorem} { Let $A$ be a $2$-atom of a reflexive locally finite
 $2$-faithful degenerate graph $\Gamma =(V,E)$ and let  $X$ be a $2$-semi-fragment not containing $A.$ Then $|A\cap X|<2,$ if one of the following conditions holds:
 \begin{itemize}
  \item[(i)] $|A|\le \nabla (X),$
  \item[(ii)] $\nabla (A)\cap \nabla (X)\neq \emptyset.$
  \end{itemize}

In particular, the intersection of two distinct  $2$-atoms of a $2$-faithful graph has a cardinality less than   $2.$

\label{faithful}} \end{theorem}
\proof
(i) follows by Theorem \ref{lem1977},(ii).
Assume that $|A\cap X|\ge 2$ and $\nabla (A)\cap \nabla (X)\neq \emptyset.$
 Since  $\kappa_{2}=\kappa _{1}$ and by Lemma \ref{degenerate}, $A$ is a $1$-fragment and $X$ is a  $1$-semi-fragment.
 By Theorem \ref{lem1977},(ii), applied with $k=1,$ $A\cap X$ is  a $1$-fragment.
 By Lemma \ref{negative}, $|\nabla (A\cap X)|\ge |\nabla (A)|\ge 2.$ Thus, $A\cap X$
 is a $2$-fragment, and hence $A\cap X=A,$ a contradiction.\qfd

\begin{lemma}
Let  $X$ and $Y$ be two $2$-atoms of
 a  reflexive locally finite
 $2$-faithful degenerate graph   $\Gamma =(V,E).$
Then
\begin{equation}
|\partial (X\cap Y)| \le |\Gamma (X)\cap \Gamma (Y)|-|X\cap Y|\le \kappa _2, \ \text{and} \label{eqk1k2}\\
\end{equation}
\begin{equation}
 |\nabla (X)\setminus \nabla  (Y)| \le |Y\setminus X|+\kappa _2-|\partial (X\cap Y)|. \label{nabla12}
\end{equation}
\label{lem19771} \end{lemma}

\proof
We use the notations of the proof of Theorem \ref{lem1977}.
By Lemma \ref{degenerate}, $X$ and $Y$ are $1$-fragments.
We shall show that
\begin{equation}\label{1223}
|R_{12}|\ge |R_{23}|.\end{equation}

This holds by  (\ref{iia}), applied with $k=1,$ if $|\nabla (X)\cap \nabla  (Y)|\ge 1.$
Suppose that $|\nabla (X)\setminus \nabla  (Y)|=0.$ We
have

$$ |R_{11}|+|R_{12}|+ |R_{13}| = |X|\le | \nabla (X)|=
 | \nabla (Y)|=
 |R_{13}|+|R_{23}|,$$

and (\ref{1223}) holds.
Thus
\begin{align*} |\partial (X\cap Y)|&\le |\Gamma (X)\cap \Gamma (Y)|-|X\cap Y|\\ &=|R_{12}|+|R_{22}|+ |R_{12}|
\\ &\le |R_{23}|+|R_{22}|+ |R_{12}|= \kappa _2 ,\end{align*}
proving (\ref{eqk1k2}).

 By Theorem  \ref{faithful}, $|X\cap Y|=1.$
Also we have,
$$|\nabla (Y)\setminus \nabla (X)|=  |R_{13}|+|R_{23}|\le |R_{13}|+|R_{12}|=|X\setminus Y|=|X|-1,
$$ proving  (\ref{nabla12}).\qfd

\section{A description of the $2$-atoms }

\begin{theorem} Let $\Gamma =(V,E)$ be a reflexive locally finite  degenerate  
and reverse degenerate graph 
such that  $\Gamma$ and $\Gamma ^-$ are $2$-faithful graphs. Then one of the following holds:

\begin{itemize}
  \item [(i)] no vertex is incident to   three pairwise  distinct $2$-atoms and incident to
  three pairwise  distinct reverse-$2$-atoms.
  \item [(ii)]   the $2$-atom has cardinality $2$ or the reverse-$2$-atom has cardinality $2.$
\end{itemize}

\label{superatoms} \end{theorem}

\proof

Let $H$ be a $2$-atom and let $K$ be a
reverse-$2$-atom. Without loss of generality we may take $|K|\ge |H|.$

Assume that (i)  does not hold.
We may choose two distinct $2$-atoms $X, Y$ incident to the same vertex $v,$ since (i) does not hold. By Theorem  \ref{faithful}, we have
 $X\cap Y=\{v\}.$

We have $\nabla (X)\not\subset \nabla (Y)$, by Lemma
\ref{negative}. Take $w\in \nabla (X)\setminus \nabla (Y).$ We have by (\ref{eqk1k2}), applied with $X$ and $Y$ permuted,

\begin{equation}\label {eqAminusF}
 |\nabla (X)\setminus \nabla (Y)|\le |Y| -1.
\end{equation}

{\bf Case }1: $L_1\cup L_2 \subset \nabla (X),$  for some  distinct reverse-$2$-atoms $L_1$ and $L_2$ with $w\in L_1\cap L_2.$
 By Theorem  \ref{faithful}, $L_1\cap L_2=\{w\}.$
 By Lemma \ref{negative},  $\nabla (Y)$ is a reverse-$2$-semi-fragment
 and $X\subset \nabla ^- (L_1)\cap \nabla ^-(L_2).$ Take an arbitrary $i\in \{1,2\}.$ Since $w\in L_i,$ we have $L_i\not\subset \nabla (Y).$
 Since  $X\subset \nabla^- (L_i),$ we have
 $v\in \nabla^- (L_i)\cap Y.$  By Theorem \ref{faithful},(ii),  $|L_i\cap \nabla (Y)|\le 1.$

 We have using  (\ref{eqAminusF}),
 $$2|Y|-3\le 2|K|-3\le |(L_1\cup L_2)\setminus \nabla (Y)|\le |\nabla (X)\setminus \nabla (Y)|\le |Y|-1,$$

 and hence $|X|=2.$  Thus (ii) holds.

{\bf Case }2:  $\nabla (X)$ contains at most one reverse-$2$-atom.

Since (ii) fails,
there exist three pairwise distinct reverse-$2$-atoms containing $w$.
 We can now assume without loss of generality that there are distinct reverse-$2$-atoms $L,M$
 with $w\in L\cap M$ and
$$ L, M \not\subset \nabla  (X).$$

By Lemma \ref{negative},  $ X \not\subset \nabla ^-  (L),$ and $ X \not\subset \nabla ^-  (M).$

We have $ |X \cap  \nabla ^-  (L)|\le 1$ and $ |X \cap  \nabla ^-(M)|\le 1,$
 by
 Theorem \ref{faithful},(i). It follows that $$|X\cap  \Gamma ^{-}(L)\cap \Gamma ^{-}(M)|\ge |X|-2.$$

By (\ref{eqk1k2}), we have
\begin{equation} \label{eq1+kappa}
 |\Gamma ^{-}(L)\cap \Gamma ^{-}(M)|\le 1+\kappa _{-2}=1+\kappa _{-1}.
\end{equation}

Since $w\in \nabla(X),$ we have $\Gamma ^-(w)\subset \Gamma ^-(\nabla(X))=V\setminus X.$ Now we have
\begin{align*} 1+\kappa _{-1}+|X|-2&\le |\Gamma ^-(w)|+|X|-2\\ &\le |\Gamma ^-(w)|+ |X \cap ( \nabla ^-  (L)\cup \nabla ^-(M))|\\
&\le  |\Gamma ^{-}(L)\cap \Gamma ^{-}(M)|\le 1+\kappa _{-1}
\end{align*} and hence $|X|=2.$\qfd

For self-reverse-graphs, this result becomes Theorem 9.3 of \cite{hiso2007}, where the hypothesis self-reverse
is omitted. The reader my suspect this, since Corollaries 9.4  and 9.6 are self-reverse.
Also the finite case  of this result is proved in \cite{hactaa}.

\section{Vertex-transitive graphs}

Let $\Gamma =(V,E)$ be a graph.
 A function $f : V  \longrightarrow V$ will be called a {\em homomorphism }if for all $x\in V$, we have
$\Gamma (f(x))=f(\Gamma (x))$. A bijective homomorphism is called an {\em automorphism}.
The graph $\Gamma$ will be called {\em vertex-transitive}  if for all $x,y\in V,$ there is an automorphism $f$
such that $y=f(x)$.
Clearly a vertex-transitive graph is regular. It is  bi-regular if $V$ is finite.
A {\em block }  of $\Gamma$ is
 a subset $B\subset V$ such that for every automorphism $f$ of $\Gamma$, either $f(B)=B$ or $f(B)\cap B=\emptyset$.

The objects defined in the previous sections (fragments, atoms and super-fragments) are
defined using the graph structure. Therefore the image  of any of these objects by a graph automorphism is an object with the same kind. This trivial observation will be used without any reference.

Recall the following result:

\begin{theorem} \cite{hast} Let $\Gamma =(V,E)$ be a reflexive locally finite $1$-separable vertex-transitive graph.
 There is a block which is either a $1$-atom or a reverse-$1$-atom.
 In particular a graph is  a Cauchy graph if and only if its block boards  and block reverse-boards have size greater  than $\delta -2.$
\label{Cauchy} \end{theorem}

\proof
By Lemma \ref{fourf}, $\Gamma$ is $1$-faithful or reverse-$1$-faithful. If $\Gamma$ is $1$-faithful, the $1$-atom is a block, by Theorem \ref{faithful}. If $\Gamma$ is reverse-$1$-faithful, the reverse $1$-atom is a bloc, by Theorem \ref{faithful}.\qfd

\begin{theorem} Let $\Gamma =(V,E)$ be a reflexive locally finite  vertex-transitive graph such that $\Gamma$ is
  degenerate and reverse-degenerate.  Then  one of the following holds:
\begin{itemize}
\item [(i)] One of the graphs $\Gamma $ and $\Gamma ^-$ is not a Cauchy graph and either the $1$-atom or
the reverse $1$-atom   is a block,
\item [(ii)] One of the graphs $\Gamma $ and $\Gamma ^-$ is non-$2$-faithful and its reverse-$2$-super-fragment  is a block,
\item [(iii)] The graphs  $\Gamma$ and $\Gamma ^-$ are $2$-faithful Cauchy graphs
and no vertex is incident to three distinct $2$-atoms and to three distinct reverse-$2$-atoms.
 \item [(iv)] The  $2$-atoms or the reverse-$2$-atom have cardinality $2$.
\end{itemize}
\label{vtvosper} \end{theorem}

\proof

By the assumptions of the theorem, we have $\kappa _2\le \delta-1.$
Consider first the case, where $\Gamma$ is not a Cauchy graph. By Theorem \ref{Cauchy}, the $1$-atom is a block  or the reverse-$1$-atom
is a block, and thus (i) holds. Similarly the result holds if $\Gamma^-$ is not a Cauchy graph. From now on, we shall assume that
the graphs $\Gamma$ and $\Gamma^-$ are Cauchy graphs.

We have  $\delta-1 \ge \kappa _2\ge \kappa _1=\delta -1.$ Similarly,  $\delta_--1 \ge \kappa _{-2}\ge \kappa _{-1}=\delta_- -1.$

Assume first that one of the graphs $\Gamma $ and $\Gamma ^-$ is non-faithful. Then its  reverse-$2$-super-fragment  is a block,
by Theorem \ref{antiatom},(ii).
Assume now that  the graphs $\Gamma $ and $\Gamma ^-$  are faithful and that (iv) does not hold.
By Theorem \ref{superatoms}, some vertex is incident to at most two distinct $2$-atoms, or to at most two
distinct reverse-$2$-atoms, and hence
(iii) holds. Since $\Gamma$ is vertex-transitive,  no vertex is incident to three distinct $2$-atoms and to three distinct reverse-$2$-atoms, and hence
(iii) holds.\qfd

Using Lemma \ref{finiteg}, we get:

\begin{corollary} \cite{hast}
Let $\Gamma =(V,E)$ be a finite reflexive degenerate    vertex-transitive graph.  Then  one of the following holds:
\begin{itemize}
\item [(i)] The graphs $\Gamma $  is not a Cauchy graph and either the $1$-atom or
the reverse $1$-atom   is a block,
\item [(ii)] One of the graphs $\Gamma $ and $\Gamma ^-$ is non-faithful and its reverse-$2$-super-fragment  is a block,
\item [(iii)] The graphs  $\Gamma$ and $\Gamma ^-$ are $2$-faithful Cauchy graphs
and no vertex is incident to three distinct $2$-atoms and to three distinct reverse-$2$-atoms.
 \item [(iv)] The  $2$-atoms or the reverse-$2$-atom have cardinality $2$.
\end{itemize}
\label{fvtvosper} \end{corollary}

\section{Cayley graphs }

Let $G$ be a group. A {\em right-$r$-progression} is a set  of the form $\{a,ra, \cdots ,r^ja\},$ for some $r\in G.$
A {\em left-$r$-progression} is a set  of the form $\{a,ar, \cdots ,ar^j\},$ for some $r\in G.$ A set will be called an {\em $r$-progression}
if it is either a {right $r$-progression} or a {left $r$-progression}.

Let  $S$  be a
subset of   $G$. The subgroup
generated by $S$ will be denoted by $\subgp{S}$.
The graph $(G,E),$ where  $ E=\{ (x,y) : x^{-1}y \
\in S \}$ is called a {\it Cayley graph}.  It will  be denoted by
$\mbox{Cay} (G,S)$.
Put $\Gamma =\mbox{Cay} (G,S)$   and  let   $F \subset G $.
Clearly
 $\Gamma (F)=FS ,$ where $FS=\{xy: x\in F \ \mbox{and}\ y\in S\}$ is the Minkowski product of $F$ by $S$.
One may check easily that left-translations are automorphisms of Cayley graphs.
In particular,  Cayley  graphs are  bi-regular and vertex-transitive.

Recall the following easy fact:

\begin{lemma} \cite{hejc3} Let $G$ be group and let $S$ be finite generating subset with $1\in S.$   For every $a\in S,$ $\subgp{S}=\subgp{Sa^{-1}}.$ Moreover  $\mbox{Cay} (G,Sa^{-1})$ and $\mbox{Cay} (G,S)$ have the same $k$-fragments. The left-translation of a $k$-atom (resp. $k$-fragment) is a $k$-atom (resp. $k$-fragment).
\label{cayleytrans} \end{lemma}
The proof follows by an easy verification. The last part can be done directly, by observing that left translations are Cayley graph automorphisms.

The next lemma allows translating intersection properties into coset covering:

\begin{lemma} \label{tcoset}
Let $A$ be a subset of a group $G.$ Put $t=|\{x^{-1}A : a\in A\}|.$
Then   $A$ is the union of $t$ right $Q$-cosets, where  $Q=\{x:xA=A\}.$
\end{lemma}

The following result generalizes a theorem of Mann \cite{manams} in the abelian case:

\begin{theorem} \cite{hast} Let $G$ be group and let $S$ be finite generating
subset with $1\in S$ and  assume that the graph  $\Gamma=\mbox{Cay} (G,S)$
is $1$-separable. Then the  $1$-atom containing
 $1$ is a subgroup or the reverse $1$-atom containing $1$ is a subgroup.
Then   $\Gamma$  a  cauchy graph if and only, for every finite subgroup $H,$
$\min (|HS|,|SH|) \le \min (|G|,|H|+|S|-1).$
\label{Cauchycayley} \end{theorem}

\proof The result follows  by combining Lemma \ref{tcoset} and Theorem \ref{Cauchy}.\qfd

We are now ready to show that either the $2$-atoms have a nice structure or the $2$-super-fragments have structure in the
degenerate case.

   \begin{theorem} \label{2ATOMCAY}
Let $G$ be group and let $S$ be finite generating subset with $1\in S$ and  put $\Gamma=\mbox{Cay} (G,S).$
Also assume that   $\Gamma$ is degenerate   and reverse-degenerate. Then there are a finite subgroup $H$
such one of the following holds.
\begin{itemize}
 \item [(i)] $H$ is a $2$-fragment or a reverse-$2$-fragment.
 \item [(ii)] $\Gamma$ and  $\Gamma ^-$ are faithful Cauchy graphs and
 there exists  an element $a,$ such that $H\cup Ha$ is a $2$-atom or a reverse-$2$-atom.
  \end{itemize}
\end{theorem}

   \proof
 By Theorem \ref{vtvosper}, one of the following conditions holds:

\begin{itemize}
\item [(1)] One of the graphs $\Gamma $ and $\Gamma ^-$ is not a Cauchy graph.
By Theorem \ref{Cauchycayley}, the $1$-atom containing $1$ is a subgroup or the $1$-atom containing $1$ is a subgroup,
and clearly (i) holds using Lemma \ref{degenerate}.

\item [(2)] One of the graphs $\Gamma $ and $\Gamma ^-$ is non-faithful and its reverse-$2$-super-fragment  is a block.
The two cases are similar  and each of them follows from the other applied to $S^{-1}.$ Consider the case where $\Gamma$ is non-faithful and take a reverse-super-fragment $K,$ with $1\in K.$  By Lemma \ref{tcoset}, $K$ is a subgroup. Now (i) holds with $H=K.$

\item [(3)] No vertex is incident to three distinct $2$-atoms and to three distinct reverse-$2$-atoms.
The two cases are similar and each of them follows from the other applied to $S^{-1}.$ Consider the case where no vertex is incident to three distinct $2$-atoms.

Let $A$ be a $2$-atom containing $1.$ It follows that the $\{a^{-1}A; a\in A\}$ consists of $2$-atoms incident to $1.$
This family contains at most two distinct subsets. By Lemma \ref{tcoset},  $A=Q\cup Qa,$ for some $a.$
 \item [(iii)] The  $2$-atoms or the reverse-$2$-atoms have cardinality $2$. The result holds in this case with $H=\{1\}.$\qfd
\end{itemize}

The next special case is enough for most of the applications:

  \begin{theorem} \label{dl}
Let  $S$ be finite generating subset of a group $G$ with $1\in S$ and $|S|<(1-\frac{1}{p}) |G|+1,$ where $p$ denotes the smallest
 cardinality of a finite non-null subgroup of $G,$ if $G$ is a torsion group, and  $p=\infty$ otherwise.

Also assume that   $\Gamma$  is degenerate and reverse-degenerate,  where $\Gamma=\mbox{Cay} (G,S)$. Then either $S$ is a progression or
  there are a finite subgroup $H$ with $|H|\ge 2$  such one of the following holds.
\begin{itemize}
 \item [(i)]  $H$ is a $2$-fragment or a reverse-$2$-fragment,
 \item [(ii)] $\Gamma$ and  $\Gamma ^-$ are faithful Cauchy graphs and
 there exists  an element $a,$ such that $H\cup Ha$ is a $2$-atom or a reverse-$2$-atom.
 \end{itemize}

  \end{theorem}

   \proof The result holds by Theorem \ref{2ATOMCAY}, unless $\Gamma$ and $\Gamma ^-$ are Cauchy graphs and the $2$-atom has size $2$ or the reverse-atom has size $2.$
   The two cases are similar and each of them follows from the other applied to $S^{-1}.$ Consider the case where a $2$-atom
   has the form $\{1,r\}.$  We have $|\{1,r\}S|=|S|+1.$
   Decompose $S=S_1\cup \cdots, S_m,$ where $S_1, \cdots, S_m$ are  right $r$-progression such that $m$ is minimal.
   In particular, $rS_i$ contains one element not contained in $S,$ for all $1\le i \le m.$
   Without loss of generality, we may assume that $|S_1|\le \cdots\le |S_m|.$ If $m=1,$ then
   $S$ is a progression and (i) holds. Assume that $m\ge 2$ and let $K$ be the subgroup generated by
   $r.$   We have $|S_2|=|K|,$ otherwise $$|S|+1=|\{1,r\}S|\ge |\{1,r\}S_1|+|\{1,r\}S_2|+|S\setminus (S_1\cup S_2)|\ge |S|+2.$$
   In particular, $K$ is a proper subgroup. Since $|S_2|= \cdots= |S_m|=|K|,$ we have also, $KS\neq G,$ otherwise $|S|\ge |G|-|K|+1
   \ge (1-\frac{1}{p}) |G|+1,$ a contradiction.
   Now we  have  $|S|-1=\kappa _2\le |KS|-|K|=|S|+|K|-|S_1|-|K|,$ and hence $|S_1|=1.$
   In particular, $\kappa _2=|KS|-|K|,$ and thus the subgroup $K$ is a $2$-fragment. Therefore (ii) holds with $H=K.$\qfd

   The last reult generalizes a result proved in the abelian case in \cite{hejc3}, and applied to the Frobenius problem in \cite{hactaa}. Our present condition $|S|<(1-\frac{1}{p}) |G|+1,$ is sharper than the condition
$|S|< |G|/2+1,$ used in \cite{hactaa}.

We have also a description of degenerate Cayley graphs.

\begin{corollary} \label{ASTERCAY}

Let $G$ be group and let $S$ be finite generating subset with $1\in S.$ The following conditions are equivalent.

 \begin{itemize}
 \item [(i)] There is a finite subset $A$ with $|A|\ge 2,$
 $$\min (|AS|,|SA|) \le \min (|G|-2,|A|+|S|-1).$$
 \item [(ii)]
   There are a finite subgroup $H$  and an element $a,$ such that
   $$\min (|H\{1,a\}S|,|S\{1,a\}H|) \le \min (|G|-2,2|H|+|S|-1).$$
  \end{itemize}

  \end{corollary}

   \proof
   Put $\mbox{Cay} (G,S).$  Clearly,  $\Gamma ^-=\mbox{Cay} (G,S^{-1}).$
   Clearly (i) implies (ii). Using Theorem \ref{2ATOMCAY}, we see easily that (ii) implies (i).

\section{Additive Combinatorics}

Recall a well known fact:

\begin{lemma}(folklore)
Let $a,b$ be elements of a group $G$ and let $H$ be a finite subgroup of $G.$
Let
 $A,B$  be  subsets of $G$
 such that  $A\subset aH$ and $B\subset Hb.$ If  $|A|+|B|>|H|,$
 then $AB=aHb$.
\label{prehistorical}
 \end{lemma}

\begin{lemma} \label{groupfrag}
Let  $S$ be finite generating
subset of a group $G$ with $1\in S.$ Assume that the graph  $\Gamma=\mbox{Cay} (G,S)$
is  degenerate and let  $H$ be a subgroup  which is a
$2$-fragment.

Then  $S^{-1}HS =S^{-1}S\cup a^{-1}Ha,$ for some $a\in S.$

\end{lemma}

\proof

Put $|HS|=k|H|$ and take a partition $S=S_1\cup \dots\cup S_k,$ where $S_i$ is the trace of $S$ on some right coset of $H.$
We shall  assume that $|S_1|\le \dots\le |S_k|.$

Observe that $k\ge 2,$ since $H$ is a proper subgroup and since $1\in S.$
By the definitions, we have  $|S|-1\ge \kappa _2(S)=|HS|-|H|.$ Thus, $2|H|-|S_1|-|S_2|\le |HS|- |S|\le |H|-1.$
Therefore $|H|+1\le |S_1|+|S_2|.$ Now for every couple $(i,j)\in [1,k] \times
 [1,k] \setminus \{(1,1)\},$ we have $$|H|+1\le |S_1|+|S_2| \le |S_i|+|S_j|,$$
  and hence by Lemma \ref{prehistorical},
$$S^{-1}S \supset S_i^{-1}S_j=S_i^{-1}HS_j.$$ Take an element $a\in S_1.$
We have $S^{-1}S\cup a^{-1}Ha=S^{-1}HS.$\qfd

{\it Proof} of {Theorem}{\ref{kneser}}
Assume that $A$ is not a progression.

Put $S=r^{-1}A,$ where $r\in A.$ Since $S\subset A^{-1}A,$ we have $\subgp{S}\subset G.$ The other inclusion
 follows since $S^{-1}S=A^{-1}A.$
Notice that $1\in S$ and that $S$ generates $G.$ Put  $\Gamma =\mbox{Cay} (G,S).$

If  $ S^{-1}S=G$ or $ SS^{-1}=G$, then (ii) holds with  $H=G$.

If  $ |S^{-1}S|= 2|G|-1$, then $ S^{-1}S=G\setminus \{u\},$ for some $u.$
Since $ S^{-1}S$ is a self-reverse set, we have $u^2=1.$ Thus (i) holds.
Similarly (i) holds, if  $ |SS^{-1}|= |G|-1$.

So we may assume that $|S|\ge 2,$ $|SS^{-1}|,|S^{-1}S|\le |G|-2.$
 By Lemma \ref{prehistorical},  $2|S|\le |G|.$
Clearly, $\Gamma$ is degenerate and reverse-degenerate.

{\bf Claim} $G$ has a subgroup which is a $2$-fragment or a reverse $2$-fragment.

Suppose the contrary. By Theorem \ref{dl},
$\Gamma$ and  $\Gamma ^-$ are faithful Cauchy graphs and
 there exists  an element $e,$ such that $H\cup He$ is a $2$-atom or a reverse-$2$-atom, where $H$ is a non-null subgroup.
 The two cases are similar  and each of them follows from the other applied with $S^{-1}$
 replacing $S.$
So we shall deal only with the case where  $H\cup He$ is a $2$-atom.

 Since $\Gamma $ is a Cauchy graph and by the assumptions, we   have $2|S|>|S^{-1}S|\ge 2|S|-1.$ Thus, $S^{-1}$ is a $2$-fragment.

Take a partition $S=S_1\cup \dots\cup S_k,$ where $S_i$ is the trace of $S$ on some right coset of $H.$
We shall  assume that $|S_1|\le \dots\le |S_k|.$ By Lemma \ref{cayleytrans}, one may take $1\in S_1.$
Assume first that $|S_1|<|H|.$ It follows that $H\cup He$ is not a subset of $S^{-1}.$
By Theorem \ref{faithful}, $(H\cup He)\cap S^{-1}=\{1\},$  contradicting Lemma \ref{katomdegree}. Thus $|S_1|=|H|,$ and hence $HS=S.$ In particular, $|HS|-|H|\le |S|-|H|,$
and $\Gamma$ would not be a Cauchy graph, a contradiction proving the claim.

{\bf Case} 1.  $G$ has a subgroup $H$ which is a $2$-fragment.

By Lemma \ref{groupfrag}, $S^{-1}HS=S^{-1}S\cup b^{-1}Hb,$ for some $b\in S=r^{-1}A.$
Therefore,  $A^{-1}rHr^{-1}A=A^{-1}A\cup b^{-1}Hb,$ for some $b\in S=r^{-1}A.$
In particular, (ii) holds.

 {\bf Case} 2.  $G$ has a subgroup which is a reverse $2$-fragment.

 By Lemma \ref{groupfrag}, $SHS^{-1}=SS^{-1}\cup bHb^{-1},$ for some $b\in S=r^{-1}A.$
 Thus, $r^{-1}AHA^{-1}r=r^{-1}AA^{-1}r\cup bHb^{-1},$
In particular, (iii) holds.\qfd

{\bf Acknowledgement}
The author would like to thank  Professors Ben Green and Terence Tao, for calling his attention to Freiman's work and Husbands dissertation.


\begin{thebibliography}{99}



\bibitem{cauchy} A.  Cauchy,  Recherches  sur  les  nombres,  {\it J.  Ecole  polytechnique}
9(1813), 99-116.


\bibitem{davenport} H. Davenport, On the addition of residue  classes, {\it J.  London  Math.
Soc.} 10(1935), 30--32.



\bibitem{diderrich} G. T. Diderrich, On Kneser's addition theorem in groups,
{\it Proc. Amer. Math. Soc. } (1973), 443-451.







\bibitem{freiman} G. Freiman, Groups and the inverse problems of additive number theory. (Russian), {\it Number-theoretic studies in the Markov spectrum and in the structural theory of set addition} (Russian), pp. 175–183. Kalinin. Gos. Univ., Moscow, 1973.



\bibitem{hcras} Y. O. Hamidoune, Sur les atomes d'un graphe orient\'e,
 {\it C. R. Acad. Sc. Paris A}  284 (1977),   1253--1256.











\bibitem{hejc3} Y. O. Hamidoune, On subsets with a small sum in abelian groups I:
The Vosper property, {\it Europ. J. of Combinatorics} 18 (1997),
541-556.

\bibitem{halgebra} Y. O. Hamidoune, An isoperimetric method in additive theory.
{\it J. Algebra} 179 (1996), no. 2, 622--630.


 \bibitem{hast} Y. O. Hamidoune,  On small subset product in a group.
Structure Theory of set-addition,  {\it Ast\'erisque}  no. 258(1999),
xiv-xv, 281--308.

 \bibitem{hactaa} Y. O. Hamidoune, {Some results in Additive number
Theory I: The critical pair Theory}, Acta Arith. 96, no. 2(2000),
97-119.




\bibitem{hiso2007} Y. O. Hamidoune, Some additive applications of the
isopermetric approach, {\it Annales de l' Institut Fourier} 58(2008),fasc. 6, 2007-2036.
 .





\bibitem{hkemp+1} Y. O. Hamidoune, A Structure Theory for Small Sum Subsets, Preprint, 2009.

\bibitem{hkft} Y. O. Hamidoune,Two Inverse results related to a question of Tao, Prprint 2010.

\bibitem{husb}L. Husbands, Approximate Groups in Additive Combinatorics: A Review of Methods and Literature,
{\it Master's dissertation}, University of Bristol, September, 2009.


\bibitem{kempacta} J. H. B. Kemperman, On small sumsets in Abelian groups,
{\it Acta Math.} 103 (1960), 66--88.


\bibitem{knesrcomp} M. Kneser, Summenmengen in lokalkompakten  abelesche Gruppen,
{\it Math. Zeit.} 66 (1956), 88--110.




\bibitem{manams} H. B. Mann,  An addition theorem for sets of elements of an Abelian group,{\it Proc. Amer. Math. Soc.} 4 (1953), 423.





\bibitem{olsonsdif} J. E. Olson, On the symmetric difference of two sets in a group,
 {\it Europ. J. Combinatorics}, (1986), 43--54.



\bibitem{sz}  O. Serra and G. Z\'emor,  Large sets with small doubling modulo $p$ are well covered by an arithmetic progression.  {\em Ann. Inst. Fourier (Grenoble)}  59  (2009),  no. 5, 2043--2060.

\bibitem{t1}T. Tao, Open question: noncommutative Freiman theorem,{\em http:{//}terrytao. wordpress.com /2007/03/02/open-question-noncommutative-freiman-theorem}.
\bibitem{t2}T. Tao, {An elementary non-commutative Freiman theorem}, {\em http://terrytao. wordpress.com /2009/11/10/an-elementary-non-commutative-freiman-theorem.}



\bibitem{tv} T. Tao, V. H. Vu,  {\it Additive Combinatorics}, Cambridge Studies
in Advanced Mathematics 105 (2006), Cambridge University Press.




\bibitem{vosper1} G. Vosper, The critical pairs of subsets of a group of prime order,
 {\it J. London Math. Soc.} 31 (1956), 200--205.


\bibitem{vw} V. H. Vu and P. M. Wood,
The inverse Erd\H{o}s-Heilbronn problem,
{\it Electron. J. Combin.} 16 (2009), no. 1, Research Paper 100, 8 pp.







\end{thebibliography}
\end{document}